\documentclass[12pt]{article}

\usepackage{latexsym}
\usepackage{epic}
\usepackage{url}
\usepackage[leqno]{amsmath}
\usepackage{amssymb}
\usepackage{mathabx}
\usepackage[all]{xy}
\usepackage{color}
\usepackage{graphicx}
\newtheorem{lemma}{Lemma}

\newcommand{\pf}{\noindent{\bf Proof: }}
\newcommand{\epf}{\hfill\hbox{\rule{6pt}{6pt}}\\}

\newcommand{\0}{{\emptyset}}

\newcommand{\p}{{\bf{p}}}
\newcommand{\Pp}{{\bf{P}}}
\newcommand{\Ptn}{{\bf{Ptn}}}

\newcommand{\Ra}{\Rightarrow}
\newcommand{\Lra}{\Longrightarrow}

\begin{document}

\title{Characterizing Block Graphs in Terms of their Vertex-Induced Partitions}

{\tiny
\author{\bf{Andreas Dress}\\
\small{CAS-MPG Partner Institute and Key Lab for Computational Biology / SIBS / CAS}\\
\small{320 Yue Yang Road, Shanghai 200031, P.\,R.\,China}\\
\small{email: andreas@picb.ac.cn}\\
\bf{Katharina T. Huber}\\
\small{University of East Anglia, School of Computing Sciences}\\
\small{Norwich, NR4 7TJ, UK}\\
\small{e-mail: Katharina.Huber@cmp.uea.ac.uk}\\ 
\bf{Jacobus Koolen}\\
\small{School of Mathematical Sciences, University of Science and Technology of China}\\
\small{96 Jinzhai Road, Hefei, Anhui 230026, P.R. China}\\
\small{email: koolen@ustc.edu.cn}\\
\bf{Vincent Moulton}\\
\small{University of East Anglia, School of Computing Sciences}\\
\small{Norwich, NR4 7TJ, UK}\\
\small{e-mail: vincent.moulton@cmp.uea.ac.uk}\\ 
\bf{Andreas Spillner}\\
\small{University of Greifswald, Department of Mathematics and Computer Science}\\
\small{17489 Greifswald, Germany}\\
\small{e-mail: andreas.spillner@uni-greifswald.de}\\
}}

\date{\today}

\maketitle

\newpage

\begin{abstract}
Given a finite connected simple graph $G=(V,E)$ with vertex set $V$ and 
edge set $E\subseteq \binom{V}{2}$, we will show that 
\vglue 1mm
\noindent
$1.$ the (necessarily unique!) smallest 
block graph with vertex set $V$ whose edge set contains $E$ is uniquely determined by the $V$-indexed family 
$\Pp_G:=\big(\pi_0(G^{(v)})\big)_{v \in V}$ of the various partitions 
$\pi_0(G^{(v)})$ of the set $V$ into the set of connected components 
of the graph $G^{(v)}:=(V,\{e\in E: v\notin e\})$, 
\vglue 1mm  
\noindent
$2.$ the edge set of this block graph coincides with set of all $2$-subsets $\{u,v\}$ of $V$ for which $u$  and $v$ are, for all $w\in V-\{u,v\}$, contained in the same connected component of $G^{(w)}$, 
\vglue 1mm
\noindent
$3.$ and an arbitrary $V$-indexed family 
$\Pp=(\p_v)_{v \in V}$ of partitions $\pi_v$ of the set $V$ is of the form $\Pp=\Pp_G$ for some connected simple graph $G=(V,E)$ with vertex set $V$ as above if and only if, for any two distinct elements $u,v\in V$, the union of the set  in $\p_v$ that contains $u$ and the set in $\p_u$ that contains $v$ coincides with the set $V$,
and $\{v\}\in \p_v$ holds for all $v \in V$. 
\vglue 1mm
\noindent
As well as being of inherent interest to the theory of block graphs,
these facts are also useful in the analysis of
compatible decompositions and block realizations of
finite metric spaces.
\end{abstract} 
{\bf Keywords and Phrases:} 
block graph, vertex-induced partition, phylogenetic combinatorics,
 compatible decomposition, strongly compatible decomposition
\section{Introduction}
\label{section:introduction}

Block graphs are a generalization of trees that arise in various areas, including metric graph theory \cite{BM86} and phylogenetics \cite{DHS97}. More specifically, a block graph is a graph in which every maximal 2-connected subgraph or block is a clique \cite{BM86,H63}. In previous work, it was shown that block graphs have interesting connections to graph-theoretical invariants of topological spaces \cite{DHKM07b} as well as to
``realisations'' of finite metric spaces in terms of weighted graphs \cite{DHKM07f}. This was part of a broader investigation into compatible decompositions of
metric spaces \cite{DHKM07d,DHKM07e,DHKM07c} which commenced in \cite[Section 4]{DMT96} based on the classical papers \cite{Isbell64,K35} and the paper \cite{BD92}. As well as forming an important part of these continued investigations 
(see e.g.\,\cite{DHKM07g,DHKM07h}) and contributing to the tasks of phylogenetic combinatorics \cite{DHKMS11}, the main result of this paper sheds light on some intriguing set-theoretical properties of block graphs and related structures.

From now on, we will consider connected simple graphs $G$ with a fixed
finite vertex set $V$. Following \cite{DHKM07f}, we will use the following notations and definitions:
\\
\\
$\bullet$  Given any set $Y$, we denote 
\\
-- by $Y\!-y$ the complement $Y\!-\{y\}$ of a one-element subset $\{y\}$ of $Y$, 
\\
-- by $\Ptn(Y)$ the set of all {\bf partitions} of $Y$  or, for short, all {\bf $Y$-partitions},
\\
--  and by $\p[y]$, for any $Y$-partition $\p$ and any element $y\in Y$, that subset $Z\in \p$ of $Y$ which contains $y$.
\\
\\
$\bullet$  Further, given a simple graph $G$ with vertex set $V$ and edge set $E\subseteq \binom{V}{2}$,
we denote 
\\
-- by $\pi_0(G)$ the $V$-partition formed by the connected components of $G$, 
\\
-- by \([G]\) the ``smallest'' block graph with vertex set \(V\) that contains \(G\) as a subgraph, i.e, the graph $(V,[E])$ with vertex set \(V\) whose edge set $[E]$ is the union of $E$ and all $2$-subsets $\{u,v\}$ of $V$ that are 
 contained in a circuit of $G$ (cf. \cite{H63}),
\\
-- by  $G[v]:=\pi_0(G)[v]$, for any vertex $v\in V$ of $G$, the connected component of $G$ containing $v$,
\\
-- by $G^{(v)}$ the largest subgraph of $G$ with vertex set $V$ for which $v$ is an isolated vertex, that is, the graph with vertex set $V$ and edge set $\{e\in E:v\notin e\}$,
\\
-- and by $\Pp_G$ the $V$-indexed family 
\begin{equation}\label{PpG}
\Pp_G:=\big(\pi_0(G^{(v)})\big)_{v \in V}\in \Ptn(V)^V
\end{equation}
of  partitions of $V$.
\\
\\
Furthermore, we define two graphs $G$ and $G'$ with vertex set $V$ to be {\bf block-equivalent} if and only if the associated block graphs $[G]$ and $[G']$ coincide. Clearly, every graph $G$ is block-equivalent to exactly one block graph, viz.\,the graph $[G]$.
  \\
 \\
And, finally, we denote by 
  \vglue 1.5mm 
\hglue 12mm   ${\bf B}(V) :=\{(V,E):E \subseteq \binom{V}{2}, (V,E) \text{ a connected block graph}\}$
  \vglue 1.5mm 
  \noindent
the set of connected block graphs with vertex set $V$\!,
 \vglue 1.5mm 
   \noindent
we define a $V$-indexed family 
  \vglue 1.5mm 
\hglue 12mm 
$\Pp=\Pp_V=(\p_v)_{v \in V}\in\Ptn(V)^V$
  \vglue 1.5mm 
  \noindent
of partitions $\p_v\in \Ptn(V)\,\,(v\in V)$ of the set $V$ to  
be a
{\bf compatible family of $V$-partitions} if 
$$ 
\p_v[u]\cup \p_u[v]=V
$$ 
holds for any two distinct elements $v, u$ in $V$, 
and $\{v\}\in \p_v$ for all $v \in V$, and we denote by 
$$
\Pp(V):=\{\Pp_V: \Pp_V \text{ is a compatible family of $V$-partitions}
\}
$$ 
the set of compatible families of $V$-partitions $\Pp_V=(\p_v)_{v \in V}$.
\\
\\ 
Note that if a $V$-indexed family of $V$-partitions is compatible,
then every pair of partitions in this family
is strongly compatible in the sense defined in \cite{DHS97}.
\\
\\
In this note, we establish the following fact:

\vglue 2mm
\noindent
{\bf Theorem:}\,\,{\em
Associating to each connected simple graph $G=(V,E)$ with vertex set $V$ the $V$-indexed family 
$\Pp_G$ as defined above, induces a one-to-one map from the set ${\bf B}(V)$ of connected block graphs with vertex set $V\,($or, equivalently, from the set of block-equivalence classes of connected simple graphs $G$ with that vertex set$)$ onto the set $\Pp(V)$ 
whose inverse is given by associating, to each family $\Pp=(\p_v)_{v \in V}$ in
$\Pp(V)$, the graph $B_\Pp:=(V,E_\Pp)$ with vertex set $V$ and edge set
$$
E_\Pp:=\big\{\{u,v\}\in \binom{V}{2}: \forall_{w\in V-\{u,v\}}\p_w[u]=\p_w[v]\big\}.
$$
In particular, given a connected graph $G=(V,E)$, the edge set $[E]$ of the associated block graph $[G]$ coincides with the set of all $2$-subsets $\{u,v\}$ of $V$ for which $G^{(w)}[u]=G^{(w)}[v]$ holds for all $w\in V-\{u,v\}$.
And given any family $\Pp=(\p_v)_{v \in V}\in \Pp(V)$, one has 
$\pi_0(B_{\Pp}^{(v)})=\p_v$ for every element $v\in V$.}

\vglue 3mm
\noindent
{\bf Proof:}
It is easy to see that, given any connected simple graph $G=(V,E)$ with vertex set $V$, the $V$-indexed family $\Pp_G=\big(\pi_0(G^{(v)})\big)_{v \in V}$ is a compatible family of $V$-partitions: Indeed, one has obviously $\pi_0(G^{(v)})[v]=\{v\}$ for every $v\in V$, and one has $\pi_0(G^{(v)})[u]\cup \pi_0(G^{(u)})[v]=V$ for any two distinct elements $v, u$ in $V$ as, given any vertex $w\in V$, there must exist a path ${\mathfrak p}=(u_0:=u,u_1, \dots, u_k:=w)$ connecting $u$ and $w$ in $G$ implying that $w\in \pi_0(G^{(v)})[u]$ holds in case $v\notin\{u_1,u_2 \dots, u_k\}$ and $w\in \pi_0(G^{(u)})[v]$ in case  $v\in\{u_1,u_2, \dots, u_k\}$.
\\
\\
We also have
 $[E]\subseteq E_{\Pp_G}$ for every connected graph $G=(V,E)$, that is, 
 $G^{(w)}[u]=G^{(w)}[v]$ holds for every
 edge $\{u,v\}\in [E]$ and all $w\in V-\{u,v\}$ because this holds clearly for every 
  edge $\{u,v\}\in E$, and it holds also for any two elements $u,v$ that are contained in a circuit 
  of $G$ as, given any vertex  $w\in V-\{u,v\}$, at least one of the two arcs of that circuit connecting $u$ and $v$ provides a path in $G^{(w)}$ connecting these two 
 vertices in that graph.
 \\
 \\
And we have  $E_{\Pp_G}\subseteq[E]$, that is, every $2$-subset $\{u,v\}$ of $V$ with $G^{(w)}[u]=G^{(w)}[v]$ for all $w\in V-\{u,v\}$
 is either an element of $E$ or contained in the vertex set of a circuit of $G$:
Indeed, employing induction relative to the length $k$ of a shortest path ${\mathfrak p}=(u_0:=u, u_1, \dots, u_k:=v)$ from $u$ to $v$ in $G$, there is nothing to prove in case $k=1$. And in case $k=2$, a circuit of $G$ containing $u$ and $v$ can be found by concatenating $p$ with a shortest path ${\mathfrak p'}=(u_0':=u, u_1', \dots, u'_{k'}:=v)$ from $u$ to $v$ in $G^{(u_1)}$ which must exist in view of $G^{(u_1)}[u]=G^{(u_1)}[v]$.
 \\
And, finally, in case $k>2$, we first observe that $G^{(w)}[u_{k-1}]=G^{(w)}[u]$ holds for all $w\in  V-\{u,u_{k-1}\}$. Indeed, in view of $\{u_{k-1},v\}\in E$, we have 
\vglue 1mm 
\hglue 20mm $G^{(w)}[u_{k-1}]=G^{(w)}[v]=G^{(w)}[u]$ 
\vglue 1mm 
\noindent
for all 
 $w\in  V-\{u,v,u_{k-1}\}$, and we have also $G^{(w)}[u]=G^{(w)}[u_{k-1}]$ for $w:=v$ in view of the fact that $(u, u_1, \dots, u_{k-1})$ is a path in $G^{(v)}$ connecting $u$ and $u_{k-1}$. 
 \\
 So, as $k>2$ implies that $\{u,u_{k-1}\}\not\in E$ must hold, our induction hypothesis implies that there must exist a circuit ${\mathfrak c}_0=(C,F)$ in $G$ with vertex set $C\subseteq V$ and edge set $F\subseteq E$ that passes through $u$ and $u_{k-1}$, i.e., with $u,u_{k-1}\in C$. Furthermore, there must exist a shortest path $(v_0:=v, v_1, \dots, v_j:=u)$ connecting $v$ and $u$ in $G^{(u_{k-1})}$. 
 \\
 Now, let $i$ denote the smallest index in $\{0,1,\dots,j\}$ with $v_i\in C$ which must exist in view of $v_j=u\in C$. In case $i=0$, we have $v=v_0\in C$ and $u\in C$ implying that $C$ is a circuit in $G$ that  passes through $u$ and $v$, as required. 
 \\
 Otherwise, we may view ${\mathfrak c}_0$ as the concatenation of 
two edge-disjoint paths, 
 \vglue 2mm
 (i) the path $\mathfrak p_0$  from $u_{k-1}$ to $v_i$ not passing 
through $u$ (unless $v_i = u$) and 
   \vglue 2mm
  (ii) the path $\mathfrak p_1$ from $v_i$ back to $u_{k-1}$ passing through $u$, 
  \vglue 2mm
 \noindent
 and then note that, replacing the path $\mathfrak p_0$ by the path $\mathfrak p_0'=(u_{k-1},v,v_1, \dots, v_i)$ (that is, concatenating $\mathfrak p_0'$ rather than $\mathfrak p_0$ with the path $\mathfrak p_1$), we obtain a new circuit $\mathfrak c_1$ in $G$ that, starting, say, in $u_{k-1}$, runs along $\mathfrak p_0'$ via $v$ over to $v_i$ and then follows the path $\mathfrak p_1$ from $v_i$ via $u$ back to $u_{k-1}$ and, thus, passes through both, $u$ and $v$, as required.
 \\
 This shows that the map from ${\bf B}(V)$ into the set $\Pp(V)$ given by associating to each connected simple graph $G=(V,E)$ with vertex set $V$ the $V$-indexed family $\Pp_G$ is a well-defined injective map, and that 
 $B_{\Pp_G}=(V,E_{\Pp_G})=(V,[E])=[G]$ holds for every connected graph $G=(V,E)$.
\vglue 2.5mm \noindent
To establish the theorem, it therefore remains to show that, conversely,
$\Pp_{B_\Pp}=\Pp$ holds for every compatible family $\Pp$ of $V$-partitions. So, assume that $\Pp$ is a fixed compatible family $\Pp=(\p_v)_{v \in V}$ of $V$-partitions. We have to show that $\p_v[u]=B_\Pp^{(v)}[u]$ holds for any two distinct elements $u,v\in V$. To this end, let us say that an element $w\in V$ {\bf separates} two elements $u,v\in V$ (relative to $\Pp$) or, for short, that ``$u|w|v$'' holds if and only if 
$w\neq u,v$ and 
$\p_w[u]\neq\p_w[v]$ (and, therefore, also $u\neq v$) holds. Clearly, one has $\{u,v\}\in E_\Pp$ for two distinct elements $u,v\in V$  if and only if there is no $w\in V-  \{u,v\}$ that separates $u$ and $v$. So, we also have 
$B_\Pp^{(v)}[u]\subseteq \p_v[u]$ for any two distinct elements $u,v\in V$ 
since, otherwise, there would exist $u',u'' \in B_\Pp^{(v)}[u]$ with 
$\{u',u''\} \in E_\Pp$, but $\p_v[u'] \neq \p_v[u'']$.
\\
\\
To  establish the converse, note that the following also holds:
\begin{lemma}\label{EPp}
Given any three distinct elements $u,v,w\in V$, the following nine assertions all are equivalent:
\begin{itemize}
\item[\bf (i)] $w\in V$ separates $u,v\in V$, i.e., $\p_w[u]\neq \p_w[v]$ or, equivalently, ``\,$u|w|v$'' holds,
\item[\bf (ii)] $\p_w[u]$ is a proper subset of $\p_v[w]$,
\item [\bf (iii)]$\p_w[u]$ is a proper subset of $\p_v[u]$,
\item [\bf (iv)]$\p_w[u]$ is a subset of $\p_v[u]$,
\item [\bf (v)] $v\notin \p_w[u]$  holds,
\item[\bf (vi)] $\p_w[v]$ is a proper subset of $\p_u[w]$,
\item[\bf (vii)] $\p_w[v]$ is a proper subset of $\p_u[v]$,
\item[\bf (viii)] $\p_w[v]$ is a subset of $\p_u[v]$,
\item[\bf (ix)]  $u\notin \p_w[v]$  holds,
\end{itemize}
and they all imply that also 
\vglue 1mm
\noindent
{\bf (x)} $w\in \p_v[u] \cap \p_u[v]$ 
\vglue 1mm
\noindent
must hold. 
\end{lemma}
{\bf Remark:}
Note that, while the last assertion follows indeed 
from the former nine, it is not equivalent to them -- as e.g.\,the binary tree with the three leaves $u,v,w$ immediately shows.
\\
\\
\pf
It is clear that,  in view of $V=\p_w[v]  \cup \p_v[w]$ and $w\not\in  \p_w[u]$,
we have
$$
\p_w[u]\neq \p_w[v] \Ra  \p_w[u]\cap\p_w[v]=\0 
\Ra \p_w[u] \subseteq V\!-\!(\p_w[v] \cup\{w\})
$$
$$
\Ra \p_w[u]  \subsetneq \p_v[w] \Ra \p_w[u]  \subseteq \p_v[w]\Ra v\notin \p_w[u] \Ra \p_w[u]\neq \p_w[v]. 
$$
So, all these assertions must be equivalent to each other, and they imply also that $u\in\p_w[u]  \subseteq  \p_v[w]$ and, hence, $\p_v[w] =\p_v[u] $ and, therefore, also $w\in \p_v[w] =\p_v[u] $ must hold.
In other words, the implications listed above yield that
$$
{\bf (i) \iff (ii) \iff (iii)  \iff  (iv) \iff  (v)}\Lra w\in \p_v[u] 
$$ 
holds.
And, switching $u$ and $v$, we also get 
$$
{\bf (i) \iff (vi) \iff (vii)  \iff  (viii)\iff (ix) \Lra} w\in \p_u[v]
$$
and, therefore, also ``$\bf (i)\Ra \bf (x)$'', as claimed.
\epf
\vglue 1.5mm
\noindent
Clearly, the lemma implies
\vglue 2.5mm
\noindent
{\bf (1)} Given any four elements $u,u',v,v'\in V$ with $u'\neq v'$ and $u\neq v$,
one has $\p_{v'}[u']\subseteq\p_v[u]$ if and only if $\p_v[u']=\p_v[u] $ and either $v=v'$ or $v|v'| u'$ holds.
\vglue 2mm
\noindent
Indeed, $\p_{v'}[u']\subseteq\p_v[u]$ implies $u' \in \p_v[u]$ as well as $v\notin \p_{v'}[u']$ and, therefore, $\p_v[u']=\p_v[u]$ as well as $v=v'$ or $v|v'| u'$ in view of ``$\bf (v) \Ra (i)$'' while, conversely, 
$\p_v[u']=\p_v[u]$ and $v=v'$ or $v|v'| u'$ implies 
$\p_{v'}[u']\subseteq\p_v[u'] = \p_v[u].$ \epf 
\vglue 2.5mm
\noindent
{\bf (2)}
Given any three {\bf distinct} elements 
$u,v,w\in V$, one has $\p_u[w]\neq\p_v[w]$.
\vglue 2mm
\noindent
Indeed, one has $\p_u[w]\neq\p_v[w]$ for any three distinct elements 
$u,v,w$ in $V$ as $\p_u[w]=\p_v[w]$ would imply $u\notin \p_v[w]$ as well as 
 $v\notin \p_u[w]$ and, therefore, $u|v|w$ as well as  
$v|u|w$ or, equivalently, $w\notin\p_v[u]$ and $w\notin\p_u[v]$
in contradiction to $V=\p_v[u]\cup\p_u[v]$.
\epf
\vglue 2.5mm
\noindent
{\bf (3)}
Next, one has $\{u,v\}\in E_\Pp$ for two distinct elements $u,v\in V$ if and only if 
$\p_v[u]$ is a minimal set in the collection 
$$
\Pp_{\p_u[v]\cap \p_v[u]}[u]:=\{\p_{w}[u]: w\in\p_u[v]\cap \p_v[u]\}
$$ 
of subsets of $V$ or, equivalently, in the collection 
$$
\Pp_{\p_u[v]}[u]:=\{\p_{w}[u]: w\in\p_u[v]\}
$$ 
or, still equivalently, in 
$$
\Pp[u]:=\{\p_{w}[u]: w\in V-u\}.
$$
Indeed, our definitions and the facts collected above imply that
\begin{eqnarray*}
\{u,v\}\not\in E_\Pp &\iff & \exists_{w\in V-\{u,v\}}\,\,\,\p_w[u] \neq\p_w[v] \,\,\,\quad\quad\quad\text{(by definition)}
\\
&\iff & \exists_{w\in \p_v[u]\cap \p_u[v]}\,\,\,\p_w[u]  \subsetneq \p_v[u] 
\quad\text{(in view of } ``{\bf (i)}\Ra {\bf (iii)}\text{''})
\\
&\iff & \p_v[u] \not\in \min\big(\Pp_{\p_u[v]\cap \p_v[u]}[u]\big)
\end{eqnarray*}
holds for any two distinct elements $u,v\in V$,
$$
\p_v[u] \not\in \min\big(\Pp_{\p_u[v]\cap \p_v[u]}[u]\big)\Lra\p_v[u] \not\in 
\min \big(\Pp_{\p_u[v]}[u]\big)\Lra\p_v[u] \not\in\min\big(\Pp[u]\big)
$$
holds for trivial reasons, and the last remaining implication
$$ 
\p_v[u] \not\in \min\big(\Pp[u]\big)\Lra \p_v[u] \not\in\min\big(\Pp_{\p_u[v]\cap \p_v[u]}[u]\big)
$$
follows from the fact that
$w\in V-u$ and $\p_w[u]  \subsetneq \p_v[u]$ implies $w\neq u,v$ as well as $u|w| v$ and, therefore, also $w\in \p_u[v] \cap \p_v[u]$ in view of ``{\bf (i)}$\Ra${\bf (x)}'', implying that also 
$$ 
\p_v[u] \not\in \min\big(\Pp_{\p_u[v]\cap \p_v[u]}[u]\big) \iff  \p_v[u] \not\in 
\min\big(\Pp_{\p_u[v]}[u]\big) $$
must hold. So,
\begin{eqnarray*}
\{u,v\}\in E_\Pp &\iff & \p_v[u]\in \min\big(\Pp_{\p_u[v]\cap \p_u[v]}[u]\big)
\\
&\iff & \p_v[u]\in \min\big(\Pp_{\p_u[v]}[u]\big)\\
&\iff & \p_v[u]\in \min\big(\Pp[u]\big) 
\end{eqnarray*}
must hold, as claimed. \epf
\vglue 2.5mm
\noindent
{\bf (4)}
Next, given three distinct elements $u,v,w\in V$ with $\{u,w\}, \{w,v\}\in E_\Pp$,
one has $\{u,v\}\in E_\Pp$ if and only if $\p_w[u]=\p_w[v]$ holds.
\vglue 2mm
\noindent
Indeed, 
 $\{u,w\}, \{w,v\}\in E_\Pp$ implies that $\p_{w'}[u]=\p_{w'}[w]=\p_{w'}[v]$ holds for all $w'\in V-\{u,v,w\}$ and that, therefore, $\{u,v\}\in E_\Pp$ or, equivalently, 
 ``$\forall_{w'\in V-\{u,v\}}\p_{w'}[u]=\p_{w'}[v]$'' holds 
if and only if one has $\p_{w'}[u]=\p_{w'}[v]$ also for the only element $w'\in V-\{u,v\}$ not in $V-\{u,v,w\}$, i.e., for $w':=w$.
\epf
\vglue 2.5mm
\noindent
{\bf (5)}
And finally, given any two distinct elements $u,v\in V$, and any sequence ${\mathfrak p}:=(u_0:=u, u_1,\dots, u_n:=v)$ of elements of $V$ such that 
$$ 
\p_{u_1}[u]\subsetneq  \p_{u_2}[u] \subsetneq \dots \subsetneq  \p_{u_n}[u]=\p_v[u]
$$ 
is a {\bf maximal} chain of subsets of $\p_v[u]$ in 
$$
\Pp^{\subseteq \p_v[u]}[u]:=\{\p_w[u]:w\in V - u, \p_w[u]\subseteq \p_v[u]\}
$$ 
ending with $\p_{v}[u]=\p_{u_n}[u]$, the sequence ${\mathfrak p}$ forms a path from $u$ to $v$ in the graph $B_\Pp=(V,E_\Pp)$, i.e., the $2$-subsets $\{u_0,u_1\},\{u_1,u_2\},\dots,\{u_{n-1},u_n\}$ of $V$ are all contained in 
$E_\Pp$. Moreover, one has $u_i | u_j | u_k$ for all $i,j,k\in \{0, 1,\dots, n\}$ with $i<j<k$ and, therefore, also $ u_1,\dots, u_{n-1}\in \p_u[v]\cap  \p_v[u] $. In particular, we must have $u|u_j|v$ for all $j \in\{1,\dots, n-1\}$ and $\p_{u_j}[u]=\p_{u_j}[u_i]$
 and $\p_{u_{i}}[v]=\p_{u_{i}}[u_{j}]$
for all $i,j=1,\dots,n$ with  $i<j$.
\vglue 2mm
\noindent
Indeed, our assumption that $\p_{u_j}[u]\subsetneq  \p_{u_k}[u]$ holds for all $j,k\in \{1, 2,\dots, n\}$ with $j<k$ implies, in view of ``$\bf (iii) \Ra \bf (i)$''  that also $u|u_j|u_k$ and, therefore, also  $\p_{u_k}[u_j]= \p_{u_k}[u]$ must hold for all $j,k=1,2,\dots,n$ with $j<k$.
 In consequence, we must also have 
 $\p_{u_j}[u_i] =\p_{u_j}[u]\subsetneq  \p_{u_k}[u]=\p_{u_k}[u_i]$ and, therefore, also  $u_i|u_j|u_k$ as well as $\p_{u_k}[u_i]= \p_{u_k}[u_j]$ for all $i,j,k\in \{0, 1,\dots, n\}$ with $i<j<k$. In particular, we must have $u|u_j|v$ for all $j \in\{1,\dots, n-1\}$ and, hence, $u_1,\dots, u_{n-1}\in \p_u[v]\cap \p_v[u]$
 and $\p_{u_j}[u]=\p_{u_j}[u_i]$
 and $\p_{u_{i}}[v]=\p_{u_{i}}[u_{j}]$
for all $i,j=1,\dots,n$ with  $i<j$, 
 as claimed.  
 \vglue 2mm
\noindent
 To establish the remaining claim that
 $\{u_0,u_1\},\{u_1,u_2\},\dots,\{u_{n-1},u_n\}\in E_\Pp$  also  holds, note first that $\p_{u_1}[u]$ is, by assumption, a minimal set in the set system 
 $\Pp^{\subseteq \p_v[u]}[u]$ 
and, therefore, also in
$\Pp[u]$ as $w \in V - u$ and $\p_w[u] \subseteq \p_{u_1}[u]$ implies $\p_w[u] \subseteq \p_v[u]$ or, equivalently,
$\p_w[u] \in \Pp^{\subseteq \p_v[u]}[u]$ and therefore, in view of the minimality of $ \p_{u_1}[u]$ in 
 $\Pp^{\subseteq \p_v[u]}[u]$, also 
$\p_{w}[u]= \p_{u_1}[u]$ or, equivalently, $w=u_1$. So,  $\{u_0,u_1\}\in E_\Pp$ must hold.
 \vglue 1mm
\noindent
 Similarly, our choice of the elements $u_0, u_1,\dots, u_n$
implies also that
\begin{equation}\label{min}
\p_{u_i}[u]\in \min\{\p_w[u]: w\in V-u \text{ and }\p_{u_{i-1}}[u]\subsetneq \p_w[u]\subseteq \p_v[u]\}
\end{equation}
must hold for all $i=2,3,\dots,n$ and, therefore, also
\begin{equation}\label{min1}
\p_{u_i}[u]=\p_{u_i}[u_{i-1}]\in \min\big(\Pp_{\p_{u_i}[u_{i-1}] \cap \p_{u_{i-1}[u_i]}}[u_{i-1}] \big)
\end{equation}
as $w\in \p_{u_i}[u_{i-1}] \cap \p_{u_{i-1}}[u_i]$ and 
$\p_w[u_{i-1}]\subsetneq \p_{u_i}[u]=\p_{u_i}[u_{i-1}]$ would imply
 $u_i\notin \p_w[u_{i-1}]$ and $w\notin \p_{u_{i-1}}[u]$ (in view of  
 $w\in \p_{u_{i-1}}[u_i]=  \p_{u_{i-1}}[v]\neq  \p_{u_{i-1}}[u]$) and, therefore, $u_{i-1}|w|u_i$ as well as $u|u_{i-1}|w$ which, in turn, would imply 
$$
\p_{u_{i-1}}[u]\subsetneq\p_w[u]=\p_w[u_{i-1}]\subsetneq \p_{u_i}[u]\subseteq \p_v[u]
$$
in contradiction to (\ref{min}).
So,  (\ref{min1}) or, equivalently,
$\{u_{i-1},u_i\}\in E_\Pp$ must hold also for all $i\in\{2,\dots,n\}$.
\epf
\\
Now, to finalize the proof of our main result, it suffices to note that, with 
$\Pp=(\p_v)_{v \in V} \in \Pp(V)$ as above, one has $ \p_v[u]\subseteq B_\Pp^{(v)}[u]$ for any two distinct elements $u,v\in V$. Yet, given any further element $u'\in \p_v[u]$, Assertion $\bf (5)$ implies that there exist two paths ${\mathfrak p}:=(u_0:=u, u_1,\dots, u_n:=v)$ and ${\mathfrak p}':=(u_0':=u', u_1',\dots, u_{n'}':=v)$ connecting $u$ and $u'$ with $v$ in $B_{\Pp}$, and Assertion $\bf (4)$ implies that also either $u_{n-1}=u_{n'-1}'$ or 
$\{u_{n-1},u_{n'-1}'\}\in E_\Pp$ holds, implying that
there exists also a path in $B_{\Pp}^{(v)}$ from $u$ to~$u'$.
\\
\\
This finishes the proof of the theorem. \epf
\\
{\bf Remark} It might also be worth noting that a compatible family of $V$-partitions $\Pp=(\p_v)_{v \in V}\in \Ptn(V)^V $ is fully encoded by the ternary relation ``$..|..|..$''$\subseteq V^3$ as $\p_v[u]$ apparently coincides, for any two distinct elements $u,v\in V$, with the set of all $w\in V-v$ for which $u|v|w$ does not hold. Consequently, one can also record the specific properties an arbitrary ternary relation ``$..|..|..$''$\subseteq V^3$ must satisfy to correspond to some 
$\Pp\in\Pp(V)$ -- a simple task that we leave as an exercise to the reader.
\subsubsection*{Acknowledgments}

A.\,Dress thanks the CAS-MPG Partner Institute for Computational Biology (PICB) in Shanghai (China), the Max Planck Institute for Mathematics in the Sciences in Leipzig (Germany), the Faculty for Mathematics at the University of Bielefeld, Germany, and the Scientific Center of the Vitroconnect-Systems GmbH
in G\"utersloh (Germany) for their continuous support.
\\
\\
J.\,Koolen thanks the 100 Talents Program of the Chinese Academy of Sciences for support.

\end{document}